\documentclass{article}

\usepackage{arxiv}

\usepackage[utf8]{inputenc} 
\usepackage[T1]{fontenc}    
\usepackage{booktabs}       
\usepackage{nicefrac}       
\usepackage{microtype}      
\usepackage{lipsum}		
\usepackage{doi} 

\usepackage{amssymb}
\usepackage{amsfonts}
\usepackage{amsmath}
\usepackage{latexsym,epsfig}
\usepackage{graphicx} 
\usepackage[usenames,dvipsnames]{color} 
\usepackage{array} 
\usepackage{color} 
\usepackage{verbatim}
\usepackage[ruled,linesnumbered]{algorithm2e}  
\usepackage{url}
\usepackage{enumerate} 
\usepackage{rotating}
\usepackage{hyperref} 
\usepackage{afterpage,lscape} 
\usepackage[misc,geometry]{ifsym} 
\usepackage{enumitem}
\setlist[enumerate]{after={\bigskip}}
\setlist[itemize]{after={\bigskip}}

\newcommand{\EndProof}{\hspace{\stretch{1}} $\Box$}

\newcommand{\J}{\mathcal{J}}

\newcommand{\R}{\mathbb{R}}
\newcommand{\N}{\mathbb{N}}

\newtheorem{definition}{Definition}
\newtheorem{example}{Example}
\newtheorem{lemma}{Lemma}
\newtheorem{proposition}{Proposition}
\newtheorem{theorem}{Theorem}
\newtheorem{corollary}{Corollary}

\title{Totally Greedy Sequences Defined by \\ Second-Order Linear Recurrences \\ With Constant Coefficients}  

\author{ \href{https://orcid.org/0000-0002-3569-3885}{\includegraphics[scale=0.06]{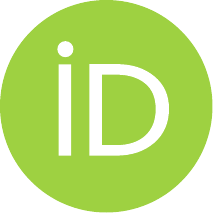}\hspace{1mm}Hebert P\'erez-Ros\'es}\thanks{Conjoint fellow, University of Newcastle, Australia} \\  \\ 
	Dept. of Computer Science and Mathematics \\       Universitat Rovira i Virgili \\ 
	Avda. Paisos Catalans 26 \\
        Tarragona 43007, Spain \\ 
	\texttt{hebert.perez@urv.cat} 
}

\renewcommand{\shorttitle}{Totally Greedy Sequences Defined by 2nd Order Linear Recurrences} 

\hypersetup{
pdftitle={Totally Greedy Sequences},
pdfsubject={ },
pdfauthor={H. P\'erez-Ros\'es},
pdfkeywords={11B37, 11B39, 11Y55, 68R05},
} 

\begin{document}
\maketitle 


\begin{abstract}
    The change-making problem consists of representing a certain amount of money with the least possible number of coins, from a given, pre-established set of denominations. The greedy algorithm works by choosing the coins of largest possible denomination first. This greedy strategy does not always produce the least number of coins, except when the set of denominations obeys certain properties. We call a set of denominations with these properties a greedy set. If the set of denominations is an infinite sequence, we call it totally greedy if every prefix subset is greedy. In this paper we investigate some totally greedy sequences arising from second-order linear recurrences with constant coefficients, as well as their subsequences, and we prove sufficient conditions under which these sequences are totally greedy. 
\end{abstract}


\section{Greedy and totally greedy sets}
\label{sec:intro}

In the \emph{money-changing problem}, or \emph{change-making problem}, we have a set of coin denominations $S = \{ s_1=1, s_2, \ldots, s_t \}$, with $s_1 < \ldots < s_t$. We also have a target amount $k$, and the goal is to make $k$ using as few coins as possible. Mathematically, we are looking for a \emph{payment vector} $(a_1, \ldots, a_t)$, such that 

\begin{align} 
	& a_i \in \N_0 \ \mbox{ for all } i = 1, \ldots, t \\
	& \sum_{i=1}^{t} a_i s_i = k, \\
	& \sum_{i=1}^{t} a_i \mbox{ is minimal,}  
\end{align} 
where $\N_0$ denotes the set of nonnegative integers. 

This problem has been extensively studied in recent years (see for instance \cite{Ada10, Cai09, Chan22, Cow08, Su23}), and it bears some relationships to other Diophantine-like problems, such as the \emph{Frobenius problem} and the \emph{postage stamp problem} \cite{Sha03}. It is also a special case of the well known \emph{knapsack problem} \cite{Ma75}. 

In regard to its computational complexity, finding the optimal payment vector for a given $k$ is NP-hard if the coins are large and represented in binary (or decimal). This result was stated by Lueker in a technical report \cite{Lue75}, and it has been cited many times, although the report itself is not available on the Internet. 

One approach for dealing with the problem is the \emph{greedy algorithm}, which proceeds by always choosing in the first place the coin of the largest possible denomination. The pseudocode of the greedy method is given in Algorithm \ref{alg:greedypayment}:  
\\\\
	\begin{algorithm}[H]
		\SetKwInOut{Input}{Input}
		\SetKwInOut{Output}{Output}
		\SetKw{DownTo}{downto}
		\vspace{.2cm}
		\Input{The set of denominations $S={1, s_2, \ldots, s_t}$, with $1 < s_2 < \ldots < s_t$, and a quantity $k \geq 0$.} 
		\Output{Payment vector $(a_1, a_2, \ldots, a_t)$.} 
		\vspace{.2cm}
		
		\For{ $i$:= $t$ \DownTo $1$} 
		{
			$a_i$ := $k$ div $s_i$\;
			$k$ :=  $k$ mod $s_i$\;
		} 
		\caption{GREEDY PAYMENT METHOD}
		\label{alg:greedypayment}  
	\end{algorithm}
\begin{definition}
	\label{def:greedycost} 
	For a given set of denominations $S={1, s_2, \ldots, s_t}$, the \emph{greedy payment vector} is the payment vector $(a_1, a_2, \ldots, a_t)$ produced by Algorithm \ref{alg:greedypayment}, and $\mbox{\textsc{GreedyCost}}_S(k) = \sum_{i=1}^{t} a_i$. 
\end{definition}

Obviously, the greedy payment vector is not necessarily optimal (i.e. $\mbox{\textsc{GreedyCost}}_S(k)$ is not always minimal among all possible payment vectors) but there exist some sets of denominations $S$ for which we can guarantee that the greedy payment vector is indeed optimal. 

\begin{example}
	\label{ex:greedyset} 
	Let $S_1 = \{ 1, 4, 6 \}$ and $S_2 = \{ 1, 2, 5 \}$ be two sets of denominations, and suppose that we want to represent the quantity $k=8$. The greedy payment vector obtained by Algorithm \ref{alg:greedypayment} with the set $S_1$ is  $(2,0,1)$, and $\mbox{\textsc{GreedyCost}}_{S_1}(8)=3$. This is obviously not optimal, since we can easily find a better way to represent the quantity $8$, namely $8=4+4$, which uses only two coins. 
	
	On the other hand, with the set $S_2$, the greedy payment vector is $(1,1,1)$, and $\mbox{\textsc{GreedyCost}}_{S_2}(8)$ is also equal to $3$. However, in this case the greedy payment vector is optimal, i.e. it is impossible to find a representation of $8$ using fewer coins of $S_2$. In fact, it can be proved that with the set $S_2$, the greedy payment vector is \emph{always} optimal for \emph{any} quantity $k$. 
\end{example}

\begin{definition}
	\label{def:greedyset} 
	If a set $S$ of denominations \emph{always} produces an optimal payment vector for \emph{any} given amount $k$, then $S$ is called \emph{orderly}, \emph{canonical}, or \emph{greedy}.
\end{definition}

Besides being interesting in their own right, greedy sets have a variety of potential applications. They can be used, for instance, to construct circulant networks with efficient routing algorithms \cite{PR21}.

There is a polynomial-time algorithm that determines whether a given set of denominations is greedy \cite{Pea05, Sha03}, as well as plenty of necessary or sufficient conditions for special families of denomination sets \cite{Ada10, Cai09, Cow08, Su23}. The current paper proceeds precisely along the latter direction. 

Obviously, a set $S$ consisting of one or two denominations is always greedy. For sets of cardinal $3$ we have the following characterization \cite{Ada10}: 

\begin{proposition}
\label{prop:triad}
The set $S = \{ 1, a, b \}$ (with $a < b$) is greedy if, and only if, $b-a$ belongs to the set 
    \begin{align*}
	 \mathfrak{D}(a) &= \{ a-1, a \} \cup \{ 2a-2, 2a-1, 2a \} \cup \ldots \{ ma-m, \ldots ma \} \cup \ldots = \\
	 &= \bigcup_{m=1}^\infty \bigcup_{s=0}^m \{ ma-s \}         
    \end{align*}      \EndProof 
\end{proposition}

The most powerful necessary and sufficient condition is given by the so-called \emph{one-point theorem} (see Theorem 2.1 \cite{Ada10}). Here we state it in a modified form: 

\begin{theorem}
	\label{theo:onepoint}
	Suppose that $S = \{ 1, s_2, \ldots, s_t \}$ is a greedy set of denominations, and $s_{t+1} > s_t$. Now let $\displaystyle m = \left\lceil \frac{s_{t+1}}{s_t} \right\rceil$. Then $\hat{S} = \{ 1, s_2, \ldots, s_t, s_{t+1} \}$ is greedy if, and only if, $\mbox{\textsc{GreedyCost}}_S(ms_t - s_{t+1}) < m$.  \EndProof 
\end{theorem}	

Notice that 
\begin{displaymath}
	(m-1)s_t +1 \leq s_{t+1} \leq ms_t, 
\end{displaymath}
by the definition of $m$. A straightforward consequence of the one-point theorem is the following

\begin{corollary}
	\label{coro:extensionbymultiple}[Lemma 7.4 of \cite{Ada10}]
	Suppose that $S = \{ 1, s_2, \ldots, s_t \}$ is a greedy set, and $s_{t+1} = us_t$, for some $u \in \N$. Then $\hat{S} = \{ 1, s_2, \ldots, s_t, s_{t+1} \}$ is also greedy.  \EndProof 
\end{corollary} 

\begin{definition}
	\label{def:totallygreedyset} 
	A set $S = \{ 1, s_2, \ldots, s_t \}$ is \emph{totally greedy}\footnote{Also called normal, or totally orderly.} if every prefix subset $\{ 1, s_2, \ldots, s_k \}$, with $k \leq t$ is greedy. 
\end{definition}

Obviously, a totally greedy set is also greedy, but the converse is not true in general. Take, for instance, the greedy set $\{ 1, 2, 5, 6, 10 \}$, whose prefix subset $\{ 1, 2, 5, 6 \}$ is not greedy. 

Nevertheless, the above situation is impossible in sets of four elements: 

\begin{proposition}[See Proposition 4.1 of \cite{Ada10}] 
    \label{prop:4elements}
    A set $S = \{ 1, s_2, s_3, s_4 \}$ is greedy if, and only if, it is totally greedy. 
\end{proposition}

Definition \ref{def:totallygreedyset} can be extended to infinite sequences in a straightforward way:   

\begin{definition}
	\label{def:totallygreedysequence} 
	Let $S = \{ s_n \}_{n=1}^{\infty}$ be an integer sequence, with $s_1 = 1$ and $s_i < s_{i+1}$ for all $i \in \N$. We say that $S$ is \emph{totally greedy} (or simply, greedy) if every prefix subset $\{ 1, s_2, \ldots, s_k \}$ is greedy.  
\end{definition}

Totally greedy sequences are mentioned very briefly in \cite{Cow08}, where some sufficient conditions are also given, that allow to construct greedy sequences from recurrence relations, although the conditions are a bit cumbersome (see Corollary 2.12 of \cite{Cow08}). In this paper we provide a simpler set of sufficient conditions that produce greedy sequences from second-order recurrences, and we also investigate some properties of these sequences and their sub-sequences. Unlike previous results, which are mainly based on combinatorial arguments, we will use a mix of analytic and combinatorial techniques that can be easily extended to other types of recurrences. We focus mainly on homogenous recurrences, though we also make some comments on the non-homogenous case.  


\section{Sequences of the form $G_{n+2} = pG_{n+1} + qG_n$} 
\label{sec:type1} 

In this section will consider sequences $\big\{ G_n \big\}_{n=1}^{\infty}$ generated by the recurrence 
\begin{equation}
	\label{eq:recurrence-type1}
	G_n = 
	\begin{cases}
		1 \mbox{ if } n=1, \\
		a \mbox{ if } n=2, \\
		pG_{n-1} + qG_{n-2}, \mbox{ if } n>2,   
	\end{cases}
\end{equation}
where $a, p, q$ are positive integers, with $a>1$, and some additional restrictions that we will see later on. 

Note that the (shifted) Fibonacci sequence $\big\{ F_n \big\}_{n=1}^{\infty}$, defined by $F_0 = 0$, $F_1 = 1$, and $F_n = F_{n-1} + F_{n-2}$, is a very special case of Equation \ref{eq:recurrence-type1}, namely, when $a = p = q = 1$. Equation \ref{eq:recurrence-type1} also generalizes the Lucas numbers (Sequence \href{https://oeis.org/A000032}{A000032} of \cite{oeis}), the Pell numbers (Sequence \href{https://oeis.org/A000129}{A000129} of \cite{oeis}), and other special cases of Horadam-type sequences.   

The parameter $a$ in Equation \ref{eq:recurrence-type1} acts as a \lq perturbation\rq \ with respect to the \lq regular\rq \ sequence 
\begin{equation}
	\label{eq:recurrence-regular-type1}
	H_n = 
	\begin{cases}
		1 \mbox{ if } n=1, \\
		p \mbox{ if } n=2, \\
		pH_{n-1} + qH_{n-2}, \mbox{ if } n>2.    
	\end{cases}
\end{equation}
Albeit small, this perturbation has a significant impact on the greediness of $\big\{ G_n \big\}_{n=1}^{\infty}$, and it also affects other properties of the sequence, as well as their proofs, which become more involved. Anyway, in addition to our natural desire for achieving greater generality, we actually \emph{need} to introduce the parameter $a$, for reasons that will become apparent in Section \ref{sec:subsequences}. 

Anyway, the characteristic polynomial associated with Equation \ref{eq:recurrence-type1} is $x^2-px-q$, with roots 
\begin{equation}
	\label{eq:roots-type1} 
	\lambda = \frac{1}{2} \left( p + \sqrt{p^2+4q} \right), \qquad \mu = \frac{1}{2} \left( p - \sqrt{p^2+4q} \right), 
\end{equation}
so that $\mu + \lambda = p$ and $\mu \lambda = -q$. 
Since the roots $\lambda$ and $\mu$ are real and distinct, the general term of $\big\{ G_n \big\}_{n=1}^{\infty}$ is 
\begin{equation}
	\label{eq:general-term-type1}
	G_{n+1} = c_1 \lambda^n + c_2 \mu^n, 
\end{equation}
where $\displaystyle c_1 = \frac{a-\mu}{\lambda-\mu}$ and $\displaystyle c_2 = \frac{\lambda-a}{\lambda-\mu}$. 

Obviously, $\big\{ G_n \big\}_{n=1}^{\infty}$ is monotonically increasing, $|\lambda| > |\mu|$, and it is also easy to verify that $\lambda > 1$. In fact, it is not difficult to prove that $\lambda > p$. At the same time, we can verify that $\mu < 0$. Now, if $q \leq p$ we can bound the roots $\lambda$ and $\mu$ with more precision. 

\begin{lemma}	
	\label{lemma:bounds-type1}
	If $\big\{ G_n \big\}_{n=1}^\infty$ is a sequence defined by Equation \ref{eq:recurrence-type1}, with $q \leq p$, and $\lambda$ and $\mu$ are the roots of the characteristic polynomial, as defined in Equation \ref{eq:roots-type1}, then 
	\begin{displaymath}
		-1 < \mu < 0 \quad \mbox{ and } \quad p < \lambda < p+1. 
	\end{displaymath}
\end{lemma} 

\textbf{Proof:} The proof is straightforward. 

\EndProof

Note that as a consequence of the above results, $c_1$ is always positive, while $c_2$ can be positive or negative, depending on $a$. In the rest of the paper, sequences that obey Equation \ref{eq:recurrence-type1}, with $q \leq p$, will also be called \emph{type-1-sequences}, and they will be the main focus of this section. 

Now, in order to apply Theorem \ref{theo:onepoint} we have to investigate the ratio 
\begin{equation}
	\label{eq:ratio-type1}
	\frac{G_{n+1}}{G_n} = \frac{c_1 \lambda^{n} + c_2 \mu^{n}}{c_1 \lambda^{n-1} + c_2 \mu^{n-1}},  
\end{equation}
where $\big\{ G_n \big\}_{n=1}^\infty$ is a type-1-sequence. 

Dividing the numerator and the denominator by $\lambda^{n-1}$ we get
\begin{equation}
	\label{eq:ratio-divided-type1}
	\frac{G_{n+1}}{G_n} = \frac{c_1 \lambda + c_2 \mu \left( \frac{\mu}{\lambda} \right)^{n-1} }{c_1 + c_2 \left( \frac{\mu}{\lambda} \right)^{n-1}}. 
\end{equation}
Since $\displaystyle \left| \frac{\mu}{\lambda} \right| < 1$, $\displaystyle \left( \frac{\mu}{\lambda} \right)^{n-1} \longrightarrow 0$, and 
\begin{equation}
	\label{eq:limit-type1}
	\lim_{n \rightarrow \infty} \frac{G_{n+1}}{G_n} = \lambda \in (p, p+1). 
\end{equation}

It will also be useful (and instructive) to investigate how the different subsequences of $\displaystyle \Big\{ \frac{G_{n+1}}{G_{n}} \Big\}$ approach the limit value of $\lambda$. 
\begin{lemma}	
	\label{lemma:monotonicity-subsequences-type1}
	Let $\big\{ G_n \big\}_{n=1}^\infty$ be a type-1-sequence. Then 
	\begin{enumerate}
	    \item If $a < \lambda$ (respectively $a > \lambda$) the subsequence $\displaystyle \Big\{ \frac{G_{2k+2}}{G_{2k+1}} \Big\}_{k=0}^\infty$ is monotonically increasing (respectively  decreasing). 
	    \item If $a < \lambda$ (respectively $a > \lambda$) the subsequence $\displaystyle \Big\{ \frac{G_{2k+1}}{G_{2k}} \Big\}_{k=1}^\infty$ is monotonically decreasing (respectively  increasing).
	\end{enumerate}
\end{lemma} 

\textbf{Proof:} One way of proving the monotonicity of the subsequence $\displaystyle \Big\{ \frac{G_{2k+2}}{G_{2k+1}} \Big\}$ is by investigating the difference 
\begin{equation}
    \label{eq:expression-first-case}
	\frac{G_{2k+2}}{G_{2k+1}} - \frac{G_{2k+4}}{G_{2k+3}} = \frac{G_{2k+2} G_{2k+3} - G_{2k+1} G_{2k+4}}{G_{2k+1} G_{2k+3}} 
\end{equation}
in the first case, and the difference 
\begin{equation}
    \label{eq:expression-second-case}
	\frac{G_{2k+1}}{G_{2k}} - \frac{G_{2k+3}}{G_{2k+2}} = \frac{G_{2k+1} G_{2k+2} - G_{2k} G_{2k+3}}{G_{2k} G_{2k+2}} 
\end{equation}
in the second case, i.e. in the subsequence $\displaystyle \Big\{ \frac{G_{2k+1}}{G_{2k}} \Big\}$. Since both denominators are positive, we will investigate the sign of the numerators 
\begin{equation}
    \label{eq:numerator-first-case}
	G_{2k+2} G_{2k+3} - G_{2k+1} G_{2k+4} = c_1 c_2 \lambda^{2k} \mu^{2k} \left( \lambda \mu^2 + \lambda^2 \mu - \mu^3 - \lambda^3 \right) 
\end{equation}
and  
\begin{equation}
    \label{eq:numerator-second-case}
	G_{2k+1} G_{2k+2} - G_{2k} G_{2k+3} = c_1 c_2 \lambda^{2k-1} \mu^{2k-1} \left( \lambda \mu^2 + \lambda^2 \mu - \mu^3 - \lambda^3 \right),  
\end{equation}
respectively. 

In the first case, the sign of the expression (\ref{eq:numerator-first-case}) depends solely on $c_2$, since $c_1$, $\lambda^{2k}$, and $\mu^{2k}$ are all positive, while $\left( \lambda \mu^2 + \lambda^2 \mu - \mu^3 - \lambda^3 \right) = -p \left( p^2+4q \right)$ is negative. If $a < \lambda$, then $c_2 > 0$, and (\ref{eq:numerator-first-case}) is negative, which means that $\displaystyle \Big\{ \frac{G_{2k+2}}{G_{2k+1}} \Big\}$ is increasing. On the other hand, if $a > \lambda$, then $c_2 < 0$, and (\ref{eq:numerator-first-case}) is positive, which means that $\displaystyle \Big\{ \frac{G_{2k+2}}{G_{2k+1}} \Big\}$ is decreasing.

In the second case, the sign of the expression (\ref{eq:numerator-second-case}) again depends solely on $c_2$, since $c_1$ and $\lambda^{2k-1}$ are positive, while $\mu^{2k-1}$ and $\left( \lambda \mu^2 + \lambda^2 \mu - \mu^3 - \lambda^3 \right)$ are negative. The rest is similar. 

\EndProof

\begin{corollary}
	\label{coro:interval1}
	Let $\big\{ G_n \big\}_{n=1}^\infty$ be a type-1-sequence. Then there exists an integer $2 \leq K_0 \leq 3$ such that for all $n \geq K_0$ we have 	
	\begin{equation}
		\label{eq:interval1}
		\frac{G_{n+1}}{G_{n}} \in (p, p+1)  
	\end{equation}
\end{corollary}

\textbf{Proof:} We just have to check that $2 \leq K_0 \leq 3$. For all $n \geq 3$ we have  
\begin{displaymath}
	\frac{G_{n+1}}{G_{n}} = \frac{pG_{n} + qG_{n-1}}{G_{n}} = p + \frac{qG_{n-1}}{pG_{n-1}+qG_{n-2}} \in (p, p+1), 
\end{displaymath}
since $q \leq p$ and $qG_{n-2} > 0$. Hence, $K_0 \leq 3$. 

Now, if additionally $a > q$, then $\displaystyle \frac{G_3}{G_2} = p + \frac{q}{a} \in (p, p+1)$, hence $K_0 = 2$.  

\EndProof 

Let's denote the prefix set $\{ 1, G_2, \ldots, G_k \}$ of $\big\{ G_n \big\}_{n=1}^\infty$ by $G^{(k)}$. We know that $G^{(2)} = \{ 1, a \}$ is always greedy, and we will now investigate when $G^{(3)}$ is greedy:  

\begin{lemma}
    \label{lemma:triad-type1}
    Let $\big\{ G_n \big\}_{n=1}^\infty$ be a type-1-sequence, then $G^{(3)} = \{ 1, a, pa+q \}$ is (totally) greedy if, and only if, $2 \leq a \leq p+q$. 
\end{lemma}

\textbf{Proof:} By Proposition \ref{prop:triad}, the set $\{ 1, a, pa+q \}$ is greedy if and only if $pa+q - a$ belongs to the set 
\begin{displaymath}
    \mathfrak{D}(a) = \{ a-1, a \} \cup \{ 2a-2, 2a-1, 2a \} \cup \ldots \{ ma-m, \ldots ma \} \cup \ldots. 
\end{displaymath}
If $a > p+q$ then $pa+q - a \notin \mathfrak{D}(a)$, so $G^{(3)}$ is not greedy. Hence $2 \leq a \leq p+q$. Let us now check that this condition is sufficient. 

We may split the condition $2 \leq a \leq p+q$ into two cases:
\begin{enumerate}
    \item $a < q$, and 
    \item $q \leq a \leq p+q$.
\end{enumerate}
In the second case it is easy to see that $pa+q - a \in \mathfrak{D}(a)$, hence $G^{(3)}$ is greedy. In the first case let $\displaystyle m' = \Big \lceil \frac{q}{a} \Big \rceil > 1$. 
\begin{align*}
    pa+q-a &= pa+q-a + (m'-1)a - (m'-1)a \\
            &= (p+m'-1)a - (m'a-q).   
\end{align*}
Thus, $pa+q - a \in \mathfrak{D}(a)$ if, and only if, $0 \leq m'a-q \leq p+m'-1$. We already know that $m'a-q \geq 0$ by the definition of $m'$. As for the other inequality, we have
\begin{displaymath}
    m'a-q < 2q-q = q \leq p < p+m'-1.  
\end{displaymath}
\EndProof 
\\\\
Now we are in the position to prove the main result of this section: 
\begin{theorem}
	\label{theo:type1}
	Let $\big\{ G_n \big\}_{n=1}^\infty$ be a type-1-sequence with $2 \leq a \leq p+q$. Then $\big\{ G_n \big\}_{n=1}^\infty$ is totally greedy.  
\end{theorem} 

\textbf{Proof:} We prove the theorem by induction. The base case is covered by Lemma \ref{lemma:triad-type1}, which guarantees that $G^{(3)}$ is greedy. So, let's suppose that $G^{(k)}$ is totally greedy for some arbitrary $k \geq 3$, and let's prove that $G^{(k+1)}$ is also greedy (and hence totally greedy). 

By Lemma \ref{lemma:bounds-type1} and Corollary \ref{coro:interval1} we know that $\displaystyle p < \frac{G_{k+1}}{G_{k}} < p+1$, so $\displaystyle m = \Big\lceil \frac{G_{k+1}}{G_{k}} \Big\rceil = p+1$. Now, 
\begin{align*}
	(p+1) G_k - G_{k+1} &= (p+1) G_k - \left( p G_k + q G_{k-1} \right) \\
	&= G_k - q G_{k-1} = \left( p G_{k-1} + q G_{k-2} \right) - q G_{k-1} \\
	&= (p-q) G_{k-1} + q G_{k-2}. 
\end{align*}
To conclude the proof, note that $\displaystyle \mbox{\textsc{GreedyCost}}_{G^{(k)}} \left( (p-q) G_{k-1} + q G_{k-2} \right) = p-q+q = p < p+1 = m$.

\EndProof 

We can now apply Theorem \ref{theo:type1} to some specific sequences: 
\begin{corollary}
	\label{coro:specific-sequences-type1}
	Consider the following sequences: 
	\begin{itemize}
		\item $\big\{ F_n \big\}_{n=1}^\infty = \{ 1, 2, 3, 5, 8, 13, \ldots \}$ (shifted Fibonacci numbers)
		\item $\big\{ P_n \big\}_{n=1}^\infty = \{ 1, 2, 5, 12, 29, 70 \ldots \}$ (shifted Pell numbers)
	\end{itemize}	
	Then, $\big\{ F_n \big\}_{n=1}^\infty$ and $\big\{ P_n \big\}_{n=1}^\infty$ are totally greedy. 
\end{corollary}

\EndProof 

If $q > p$ we can no longer guarantee that $\big\{ G_n \big\}_{n=1}^\infty$ is totally greedy. Take for instance the (shifted) Jacobstahl numbers:  $\big\{ \J_n \big\}_{n=2}^\infty = \{ 1, 3, 5, 11, 21, 43, 85 \ldots \}$, defined by $\J_0 = 0$, $\J_1 = 1$, and $\J_n = \J_{n-1} + 2 \J_{n-2}$. Indeed, $\displaystyle \Big\lceil \frac{\J_{7}}{\J_{6}} \Big\rceil = 3 = m$, but $3 \cdot 21 - 43 = 20$, and $\displaystyle \mbox{\textsc{GreedyCost}}_{\J^{(6)}} \left( 20 \right) = 4 > m$. Hence, $\J^{(7)} = \{ 1, 3, 5, 11, 21, 43 \}$ is not greedy, and $\big\{ \J_n \big\}_{n=2}^\infty$ is not totally greedy. Intuitively, the problem here seems to be the largest root $\lambda = 2$, which is an integer, while the ratio $\displaystyle \frac{\J_{k+1}}{\J_{k}}$ approaches $\lambda$ from above and below, and hence sometimes $\displaystyle \Big\lceil \frac{\J_{k+1}}{\J_{k}} \Big\rceil = 3$. 





\subsection{The non-homogenous case} 
\label{sec:non-homogenous}

Let us now consider sequences $\big\{ T_n \big\}_{n=1}^{\infty}$ defined by the non-homogenous recurrence relation  
\begin{equation}
	\label{eq:non-homogenous-type1}
	T_n = 
	\begin{cases}
		1 \mbox{ if } n=1, \\
		a \mbox{ if } n=2, \\
		pT_{n-1} + qT_{n-2} \pm r, \mbox{ if } n>2, 
	\end{cases}
\end{equation}
where $a, p, q, r$ are positive integers. 

By subtracting $T_n = pT_{n-1} + qT_{n-2} \pm r$ from $T_{n+1} = pT_{n} + qT_{n-1} \pm r$, we get the homogenous third-order recurrence $T_{n+1} = (p+1)T_{n} + (q-p)T_{n-1} - qT_{n-2}$, with characteristic equation $x^3-(p+1)x^2+(p-q)x+q=0$. The roots of this characteristic equation are the same $\lambda$ and $\mu$ of Equation \ref{eq:roots-type1}, plus the additional root $\xi = 1$. Hence, the general term of the sequence is
\begin{equation}
	\label{eq:general-term-non-homogenous}
	T_n = c_1 \lambda^n + c_2 \mu^n + c_3, 
\end{equation}
for some constants $c_1, c_2, c_3$.  

By proceeding as in the homogenous case, we get that 
\begin{equation}
	\label{eq:limit-non-homogenous}
	\lim_{n \rightarrow \infty} \frac{T_{n+1}}{T_n} = \lambda,
\end{equation}
which lies in the interval $(p, p+1)$, provided that $q \leq p$. Therefore, we get the following result, as a weaker version of Theorem \ref{theo:type1}:

\begin{proposition}
	\label{prop:non-homogenous}
	Let $\big\{ T_n \big\}_{n=1}^\infty$ be a sequence defined by Equation \ref{eq:non-homogenous-type1}, with $q \leq p$. Let $K_0$ be the integer described in Corollary \ref{coro:interval1}, and suppose that $K_0 \geq 3$, and the set $\{ 1, T_2, \ldots, T_{K_0} \}$ is greedy. Then $\big\{ T_n \big\}_{n=1}^\infty$ is totally greedy.  
\end{proposition}  

\textbf{Proof:}
The proof is similar to that of Theorem \ref{theo:type1}. 

\EndProof

In order to strengthen Proposition \ref{prop:non-homogenous} we would have to make a detailed analysis of the sequence $\big\{ T_n \big\}_{n=1}^\infty$, as we did for type-1 sequences. However, this will be deferred for future work.

\section{Sequences of the form $J_{n+2} = pJ_{n+1} - qJ_n$} 
\label{sec:type2} 

We will now consider sequences $\big\{ J_n \big\}_{n=1}^{\infty}$ generated by the recurrence 
\begin{equation}
	\label{eq:recurrence-type2}
	J_n = 
		\begin{cases}
			1 \mbox{ if } n=1, \\
			a \mbox{ if } n=2, \\
			pJ_{n-1} - qJ_{n-2}, \mbox{ if } n>2, 
		\end{cases}
\end{equation}
where $a, p, q$ are integers, with $p, q > 0$ and $a > 1$. 

Recurrences defined by Equation \ref{eq:recurrence-type2} are intimately related with the type-1 sequences that we have just seen in Section \ref{sec:type1}. This phenomenon was already noticed in the case of Fibonacci numbers in \cite{Raj04, Sil06}, and we will take a deeper look at it in Section \ref{sec:subsequences}. A special case of Equation \ref{eq:recurrence-type2} is investigated in \cite{Cri07} under the name \emph{d-sequences}. 

The characteristic polynomial associated with this recurrence is $x^2-px+q$, with roots 
\begin{equation}
	\label{eq:roots-type2} 
	\lambda = \frac{1}{2} \left( p + \sqrt{p^2-4q} \right), \qquad \mu = \frac{1}{2} \left( p - \sqrt{p^2-4q} \right), 
\end{equation}
so that $p = \mu + \lambda$ and $q = \mu \lambda$. 

If $\lambda$ and $\mu$ are real and distinct, then the general term of $\big\{ J_n \big\}_{n=1}^{\infty}$ is 
\begin{equation}
	\label{eq:general-term-type2}
	J_n = c_1 \lambda^{n-1} + c_2 \mu^{n-1}, 
\end{equation}
where $\displaystyle c_1 = \frac{a-\mu}{\lambda-\mu}$ and $\displaystyle c_2 = \frac{\lambda-a}{\lambda-\mu}$. 

The roots (\ref{eq:roots-type2}) are distinct and real if, and only if $\displaystyle q < \frac{p^2}{4}$. 

As in Section \ref{sec:type1}, we are interested in a subclass of the sequences defined by Equation \ref{eq:recurrence-type2}, namely those with $1 \leq q < p-1$, which henceforth will be called \emph{type-2-sequences}. The condition $1 \leq q \leq p-2$ implies $\displaystyle q < \frac{p^2}{4}$, and it also implies that the dominant root $\lambda > 1$, so that the sequence $\big\{ J_n \big\}_{n=1}^{\infty}$ is monotonically increasing. Additionally, we can see that $0 < \mu <1$, which means that $\lambda$ cannot be an integer. 

In fact, we can bound the root $\lambda$ even further: 
\begin{lemma}	
	\label{lemma:bounds-type2}
	If $\big\{ J_n \big\}_{n=1}^\infty$ is a type-2-sequence, and $\lambda$ and $\mu$ are the roots of the characteristic polynomial, as defined in Equation \ref{eq:roots-type2}, then $p-1 < \lambda < p$. 
\end{lemma} 

\textbf{Proof:} The proof is straightforward and is omitted here. 

\EndProof


As before, in order to apply Theorem \ref{theo:onepoint} we have to investigate the ratio 
\begin{equation}
	\label{eq:ratio-type2}
	\frac{J_{n+1}}{J_n} = \frac{c_1 \lambda^{n} + c_2 \mu^{n}}{c_1 \lambda^{n-1} + c_2 \mu^{n-1}} 
\end{equation}

Dividing the numerator and the denominator by $\lambda^n$, as in Equation \ref{eq:ratio-divided-type1} we get
\begin{equation}
	\label{eq:limit-type2}
	\lim_{n \rightarrow \infty} \frac{J_{n+1}}{J_n} = \lambda \in (p-1, p). 
\end{equation}

In type-2 sequences we are also interested in the behaviour of the subsequences of $\displaystyle \Big\{ \frac{J_{n+1}}{J_{n}} \Big\}$, and how they approach the limit value $\lambda$. 
\begin{lemma}	
	\label{lemma:monotonicity-subsequences-type2}
	Let $\big\{ J_n \big\}_{n=1}^\infty$ be a type-2-sequence. If $a < \lambda$ (respectively $a > \lambda$) then the sequence $\displaystyle \Big\{ \frac{J_{n+1}}{J_{n}} \Big\}_{n=1}^\infty$ is  monotonically increasing (respectively  decreasing). 
\end{lemma} 

\textbf{Proof:} As in Lemma  \ref{lemma:monotonicity-subsequences-type1}, we will investigate the difference 
\begin{equation}
    \label{eq:difference-of-quotients-type2}
    \frac{J_{n+1}}{J_{n}} - \frac{J_{n+2}}{J_{n+1}} = \frac{J_{n+1}^2 - J_{n} J_{n+2}}{J_{n}J_{n+1}}. 
\end{equation}
Since the denominator $J_{n}J_{n+1} > 0$, the sign of Equation \ref{eq:difference-of-quotients-type2} depends on the numerator
\begin{equation}
    \label{eq:numerator-type2}
    J_{n+1}^2 - J_{n} J_{n+2} = c_1 c_2 \lambda^{n-1} \mu^{n-1} \left( 2 \lambda \mu - \lambda^2 - \mu^2 \right).  
\end{equation}
The factors $c_1, \lambda^{n-1}, \mu^{n-1}$ are positive, while $2 \lambda \mu - \lambda^2 - \mu^2 = 4q-p^2$ is negative. Hence, the sign of (\ref{eq:numerator-type2}) depends solely on $\displaystyle c_2 = \frac{\lambda-a}{\lambda-\mu}$. The rest is straightforward.

\EndProof 

From Lemmas \ref{lemma:bounds-type2} and  \ref{lemma:monotonicity-subsequences-type2} we get 

\begin{lemma}
	\label{lemma:interval2}
	Let $\big\{ J_n \big\}_{n=1}^\infty$ be a type-2 sequence. Then there exists an integer $1 < K_0 \leq 4$ such that for all $n \geq K_0$ we have 	
	\begin{equation}
		\label{eq:interval2}
		\frac{J_{n+1}}{J_{n}} \in (p-1, p).  
	\end{equation}
	In particular, if $\displaystyle a > \frac{q}{2}$ then $K_0 \leq 3$, and if $a > q$, then $K_0 = 2$. 
\end{lemma}

\textbf{Proof:} The existence of $K_0$ follows directly from Lemmas \ref{lemma:bounds-type2} and  \ref{lemma:monotonicity-subsequences-type2}, so we only have to verify the bounds. 

In the case $a>q$ it is easy to verify that  $\displaystyle \frac{J_3}{J_2} \in (p-1, p)$. In the case $\displaystyle a > \frac{q}{2}$ we have to check that $\displaystyle \frac{J_4}{J_3} \in (p-1, p)$, which is equivalent to $\displaystyle \frac{qa}{pa-q} < 1$, which in turn amounts to $pa-q-qa > 0$. But since $p \geq q+2$ we get $pa-q-qa \geq 2a-q > 0$, as desired. 

Consequently, we are left with the case $\displaystyle 2 \leq a \leq \frac{q}{2}$, $4 \leq q \leq p-2$, whence $p \geq 6$. In order to prove that $K_0 \leq 4$ in this case, we have to check that $\displaystyle \frac{J_5}{J_4} = p-q\frac{J_3}{J_4} \in (p-1, p)$, which is tantamount to proving that $\displaystyle q\frac{J_3}{J_4} < 1$, which in turn amounts to $p^2a - pq - qa - pqa + q^2 > 0$. 

Now, to prove this inequality we can use standard techniques from multivariable Calculus. Let us fix $p$ and investigate the (continuous) bivariate polynomial $F_p(q,a) = p^2a - pq - qa - pqa + q^2$ in the triangular region $\displaystyle T = \big\{ (q,a) \in \R^2: 4 \leq q \leq p-2, 2 \leq a \leq \frac{q}{2} \big\}$. Since $T$ is closed and bounded, and $F_p(q,a)$ is continuous everywhere in $\R^2$, then by the Weierstrass extreme value theorem $F_p(q,a)$ achieves its maximum and minimum values in $T$. Then we only have to check that the minimum value is positive. 

Let us first find the stationary points of $F_p(q,a)$ in $\R^2$. We compute the partial derivatives and equate them to zero:  
\begin{align*}
    \frac{\partial F_p}{\partial q} &= -p-a-pa+2q = 0\\
    \frac{\partial F_p}{\partial a} &= p^2-q-pq = 0. 
\end{align*}
By solving the first equation for $a$ we get $\displaystyle a = \frac{2q-p}{p+1}$ (note that $p \neq -1$). By solving the second equation for $q$ we get $\displaystyle q = \frac{p^2}{p+1}$. We can then substitute $\displaystyle \frac{p^2}{p+1}$ for $q$ in the expression of $a$ to get $\displaystyle a = \frac{p(p-1)}{(p+1)^2}$. 

Hence, the only stationary point is $\displaystyle \left( \frac{p^2}{p+1}, \frac{p(p-1)}{(p+1)^2} \right)$. However, it is not difficult to verify that this point does not belong to $T$, since $\displaystyle \frac{p^2}{p+1} < p-2$ implies $p \leq -2$. 

So, let us now investigate the boundary of $T$ (including its vertices as a special case). To begin with we can substitute $\displaystyle \frac{q}{2}$ for $a$ in the expression of $F_p(q,a)$, and we get the function 
\begin{displaymath}
    \hat{F}_{p,1}(q) = \frac{1}{2} (1-p)q^2 + \left( \frac{1}{2} p^2 - p \right) q, 
\end{displaymath}
which only depends on $q$. Solving 
\begin{displaymath}
    \hat{F}'_{p,1}(q) = (1-p)q + \frac{1}{2} p^2 - p = 0,
\end{displaymath}
we get $\displaystyle q = \frac{\frac{1}{2} p^2 - p}{p-1}$. Instead of trying to find out if this point lies in $T$, we can check that $\hat{F}''_{p,1}(q) = (1-p) < 0$, which means that this point is a local maximum, and therefore irrelevant for our purposes. 

Next we check the side $q=p-2$. After the substitution we get the function 
\begin{displaymath}
    \hat{F}_{p,2}(a) = pa - 2p + 2a + 4, 
\end{displaymath} 
whose derivative is $\hat{F}'_{p,2}(a) = p+2$, which means that $\hat{F}_{p,2}(a)$ does not have any stationary points. 

On the side $a=2$ we get the function 
\begin{displaymath}
    \hat{F}_{p,3}(q) = 2p^2+q^2-3pq-2q. 
\end{displaymath}
Solving 
\begin{displaymath}
    \hat{F}'_{p,3}(q) = 2q-3p-2 = 0, 
\end{displaymath}
we get $q = \displaystyle \frac{3}{2} p - 1$, which also falls without $T$. 

Thus, we are only left with the vertices of $T$ as potential minima. Namely, 
\begin{align*}
    & F_p(p-2,2) = 8, \\
    & F_p\left( p-2,\frac{p-2}{2} \right) = \frac{1}{2}(p-2)^2 = \frac{1}{2}q^2 > 0, \\
    & F_p(4,2) = p^2-6p+4 = \left( p-(3-\sqrt{5}) \right)\left( p-(3+\sqrt{5}) \right) > 0 \mbox{ for } p \geq 6. 
\end{align*}

Consequently, $\displaystyle \frac{J_5}{J_4} \in (p-1, p)$ for $p \geq 6$, as desired. 

\EndProof 


The next lemma establishes the bounds for $\displaystyle \frac{J_4}{J_3}$ in the worst case of Lemma \ref{lemma:interval2}, i.e. when $K_0 = 4$, so that $\displaystyle \frac{J_5}{J_4} \in (p-1, p)$ but $\displaystyle \frac{J_4}{J_3} \notin (p-1, p)$. 

\begin{lemma}
    \label{lemma:interval2B}
    Let $\big\{ J_n \big\}_{n=1}^\infty$ be a type-2 sequence. Suppose also that the number $K_0$ described in Lemma \ref{lemma:interval2} is equal to $4$. Then, $\displaystyle \frac{J_4}{J_3} \in (p-2, p-1]$. 
\end{lemma}

\textbf{Proof:}  Recall from Lemma \ref{lemma:interval2} that, if $\displaystyle a > \frac{q}{2}$ then $K_0 \leq 3$, hence the condition $K_0 = 4$ implies $\displaystyle a \leq \frac{q}{2}$. That means our sequence $\displaystyle \Big\{ \frac{J_{n+1}}{J_{n}} \Big\}$ is monotonically increasing, by Lemma \ref{lemma:monotonicity-subsequences-type2}. 

Thus, we already know that $\displaystyle \frac{J_4}{J_3} < p$ but $\displaystyle \frac{J_4}{J_3} \notin (p-1, p)$, so we only have to prove that $\displaystyle \frac{J_4}{J_3} > p-2$. Now,
\begin{displaymath}
    \frac{J_4}{J_3} = p - q \frac{J_2}{J_3} = p - \frac{qa}{pa-q}, 
\end{displaymath}
and 
\begin{align*}
    p - \frac{qa}{pa-q} > p-2 \quad &\mbox{ iff } \quad \frac{qa}{pa-q} < 2 \\ 
            &\mbox{ iff } \quad 2pa-2q-pa > 0. 
\end{align*}
For proving the latter inequality we can use the same technique that we used in the proof of Lemma \ref{lemma:interval2}, i.e. we fix $p$ and we define the two-variable function $F_p(q,a) = 2pa-2q-pa$. Then we can look for the minimum value of $F_p(q,a)$ over the triangular region $\displaystyle T = \big\{ (q,a) \in \R^2: 4 \leq q \leq p-2, 2 \leq a \leq \frac{q}{2} \big\}$. Finally we just have to check that this minimum is positive. 

So, let us first find the stationary points of $F_p(q,a)$ in $\R^2$. We compute the partial derivatives and equate them to zero:  
\begin{align*}
    \frac{\partial F_p}{\partial q} &= -2-a = 0\\
    \frac{\partial F_p}{\partial a} &= 2p-q = 0. 
\end{align*}
The only solution of this system is the point $(2p,-2)$, which is not in $T$. Hence, we just need to look at the boundary of $T$ (including its vertices as a special case). 

To begin with we can substitute $\displaystyle \frac{q}{2}$ for $a$ in the expression of $F_p(q,a)$, and we get the function 
\begin{displaymath}
    \hat{F}_{p,1}(q) = pq - 2q - \frac{1}{2} q^2, 
\end{displaymath}
which only depends on $q$. Then we set  
\begin{displaymath}
    \hat{F}'_{p,1}(q) = p-2-q = 0. 
\end{displaymath}
It turns out that the solution to the previous equation, the point $\displaystyle \left( p-2, \frac{p-2}{2} \right)$, is a vertex of $T$; we will check all the vertices at the end. 
Next we check the side $q=p-2$. After the substitution we get the function 
\begin{displaymath}
    \hat{F}_{p,2}(a) = 2pa - 2p - ap + 2a + 4, 
\end{displaymath} 
whose derivative is $\hat{F}'_{p,2}(a) = p+2$, which can only be zero if $p=-2$. 

On the side $a=2$ we get the function 
\begin{displaymath}
    \hat{F}_{p,3}(q) = 4p-4q,  
\end{displaymath}
with constant derivative $\hat{F}'_{p,3}(q) = -4$, which means that there are no stationary points on this side either.  

Thus, we are only left with the vertices of $T$ as potential minima. Namely, 
\begin{align*}
    & F_p(p-2,2) = 8, \\
    & F_p\left( p-2,\frac{p-2}{2}\right) = \frac{1}{2}p^2 - 2p + 2 = \frac{1}{2} (p-2)^2 = \frac{1}{2}q^2 > 0, \\
    & F_p(4,2) = 4p-16 > 0 \mbox{ for } p > 4. 
\end{align*}

Consequently, $\displaystyle \frac{J_4}{J_3} \in (p-2, p-1]$  as desired. 

\EndProof 

The approach we have used in proving the inequalities of Lemmas \ref{lemma:interval2} and \ref{lemma:interval2B} might not be as elegant as an ad hoc argument, and is probably not the shortest route to the proof, but on the other hand, it is easy to extend to other inequalities, and exhibits nicely the interplay between Analysis and Combinatorics. 

Next we continue the study of the worst-case sequences, and in particular, we turn our attention to the prefix set $\{ 1, J_2, J_3, J_4 \}$, whose quotients $\displaystyle \frac{J_{k+1}}{J_{k}}$ lie outside $(p-1, p)$. 

As in Section \ref{sec:type1}, let's denote the prefix set $\{ 1, J_2, \ldots, J_k \}$ by $J^{(k)}$. 

\begin{lemma}
    \label{lemma:prefix2B}
    Let $\big\{ J_n \big\}_{n=1}^\infty$ be a type-2 sequence. Then the prefix set $J^{(4)}$ is (totally) greedy. 
\end{lemma}

\textbf{Proof:}  Let us first check that $J^{(3)}$ is greedy. By Proposition \ref{prop:triad}, the set  $J^{(3)} = \{ 1, a, pa-q \}$ is greedy iff the difference $pa-q-a = (p-1)a-q$ belongs to the set
\begin{displaymath}
    \mathfrak{D}(a) = \{ a-1, a \} \cup \{ 2a-2, 2a-1, 2a \} \cup \ldots \{ ma-m, \ldots ma \} \cup \ldots,   
\end{displaymath} 
i.e. iff $0 \leq q \leq p-1$, but $0 < q \leq p-2$ by the assumptions. 

By Lemma \ref{lemma:interval2}, either $K_0 \leq 3$ or $K_0 = 4$. If $K_0 \leq 3$, then $\displaystyle m = \Big\lceil \frac{J_{4}}{J_{3}} \Big\rceil = p$. Therefore,  
\begin{displaymath}
	p J_3 - J_{4} = p J_3 - \left( p J_3 - q J_{2} \right) = q J_{2},   
\end{displaymath}
and $\displaystyle \mbox{\textsc{GreedyCost}}_{J^{(3)}} \left( q J_{2} \right) = q < p-1 < m$. 

Now, in the case $K_0 = 4$ we know that $\displaystyle a \leq \frac{q}{2}$. We also know, by Lemma \ref{lemma:interval2B}, that $\displaystyle  \frac{J_4}{J_3} \in (p-2, p-1]$, hence $\displaystyle \Big \lceil \frac{J_4}{J_3} \Big \rceil = p-1$. It is trivial to check that the condition $\displaystyle \frac{J_4}{J_3} \leq p-1$ is equivalent to $qa - pa + q \geq 0$, which in turn implies $\displaystyle a \leq \frac{q}{p-q} \leq \frac{q}{2}$. 

Coincidentally, in order to prove that $J^{(4)}$ is greedy we have to compute 
$\displaystyle \mbox{\textsc{GreedyCost}}_{J^{(3)}} \left( (p-1) J_3 - J_{4} \right) = \mbox{\textsc{GreedyCost}}_{J^{(3)}} \left( qa - pa + q \right)$. There are two possible cases for $qa - pa + q$: 
\begin{align*}
    0 &\leq qa - pa + q < J_2 = a \\
    a &\leq qa - pa + q < J_3 = pa-q, 
\end{align*}
which translate to 
\begin{align*}
    \frac{q}{p-q} \geq \; & a > \frac{q}{p-q+1} \\
    \frac{q}{p-q+1} \geq \; & a > \frac{2q}{2p-q},  
\end{align*}
respectively. Note that a hypothetical third case, $pa-q \leq qa - pa + q < J_4 = p^2a-pq-qa$, is impossible, for it would imply $\displaystyle \frac{2q}{2p-q} \geq 2$, which would mean that $q \geq p$. 

So, in the first case $\displaystyle  \mbox{\textsc{GreedyCost}}_{J^{(3)}} \left( qa - pa + q \right) = qa - pa + q < a \leq \frac{q}{p-q} < p-1$. 

In the second case $qa-pa+q = qa-pa+q'a+r = (q-p+q')a + r$, where $q'$ is the quotient of the integer division of $q$ by $a$, and $0 \leq r < a$ is the remainder. So, $\displaystyle  \mbox{\textsc{GreedyCost}}_{J^{(3)}} \left( (q-p+q')a + r \right) = q-p+q'+r$. Now, recall that $a \geq 2$, hence  $\displaystyle q' \leq \frac{q}{2}$. Additionally, $\displaystyle r < a \leq \frac{q}{2}$ (in fact, $\displaystyle a \leq \frac{q}{p-q+1} < \frac{q}{2}$) and $q \leq p-2$. Combining these inequalities we get $q-p+q'+r < p-4$.  

\EndProof 

Now we are ready to prove the main result of this section: 

\begin{theorem}
	\label{theo:type2}
	Let $\big\{ J_n \big\}_{n=1}^\infty$ be a type-2 sequence, i.e. a sequence defined by Equation \ref{eq:recurrence-type2}, with $q < p-1$. Then $\big\{ J_n \big\}_{n=1}^\infty$ is totally greedy.  
\end{theorem}

\textbf{Proof:} We prove the theorem by induction. The base case is guaranteed by Lemma \ref{lemma:prefix2B}, which states that $J^{(4)}$ is (totally) greedy. So, let's suppose that $J^{(k)}$ is totally greedy for some arbitrary $k \geq 4$, and let's prove that $J^{(k+1)}$ is also greedy (and hence totally greedy). 

By Lemmas \ref{lemma:bounds-type2} -- \ref{lemma:interval2} we know that $\displaystyle p-1 < \frac{J_{k+1}}{J_{k}} < p$, so $\displaystyle m = \Big\lceil \frac{J_{k+1}}{J_{k}} \Big\rceil = p$. Now, 
\begin{align*}
	p J_k - J_{k+1} &= p J_k - \left( p J_k - q J_{k-1} \right) \\
	&= q J_{k-1}.  
\end{align*}
Finally, $\displaystyle \mbox{\textsc{GreedyCost}}_{J^{(k)}} \left( q J_{k-1} \right) = q < p-1 < m$.

\EndProof 





In summary, we have verified that type-2 sequences are totally greedy, regardless of $a$, although this fact is harder to prove than the corresponding property of type-1 sequences. In any case, this property of type-2 sequences will turn out to be very useful in Section \ref{sec:subsequences}, where we analyze the subsequences of both type-1 and type-2 sequences. 

Unfortunately, Theorem \ref{theo:type2} does not seem to have a straightforward generalization to the non-homogenous case, as Theorem \ref{theo:type1} does. Take for instance the non-homogenous sequence $\big\{ T_n \big\}_{n=1}^{\infty}$ defined by the equation  
\begin{equation}
	\label{eq:counterexample-nonhomogeneous-type2}
	T_n = 
	\begin{cases}
		1 \mbox{ if } n=1, \\
		3 \mbox{ if } n=2, \\
		3T_{n-1} - T_{n-2} + 2 \mbox{ if } n>2.  
	\end{cases}
\end{equation}
First note that the parameters $p=3$ and $q=1$ meet the condition $q<p-1$ of Theorem \ref{theo:type2}. Note also that $\lambda \approx 2.618$, and for all $n \geq 3$, the ratio $\frac{T_{n+1}}{T_{n}}$ lies in the interval $(2, 3)$, hence $K_0 = 3$. Moreover, $\{ 1, T_2, T_3, T_4 \}$ is totally greedy. However, neither $\{ 1, T_2, \ldots, T_5 \}$ nor $\{ 1, T_2, \ldots, T_6 \}$ are greedy. 

Determining the exact conditions under which Theorem \ref{theo:type2} can be generalized to the non-homogenous case remains an open problem. 

\section{Even and odd subsequences}
\label{sec:subsequences}

Given a totally greedy type-1-sequence $\big\{ G_n \big\}_{n=1}^{\infty}$, or a totally greedy type-2-sequence $\big\{ J_n \big\}_{n=1}^{\infty}$, we now want to study the behaviour of the subsequence formed by the even terms, i.e. $\big\{ G_{2k} \big\}_{k=1}^{\infty}$ and $\big\{ J_{2k} \big\}_{k=1}^{\infty}$, and the subsequence formed by the odd terms, i.e. $\big\{ G_{2k-1} \big\}_{k=1}^{\infty}$ and $\big\{ J_{2k-1} \big\}_{k=1}^{\infty}$. For simplicity we may call them the \emph{even} and \emph{odd subsequences}, respectively. 

Note in passing that starting the sequences at $G_1 = J_1 = 1$, rather than starting at $G_0 = J_0 = 1$, was a convenient choice from the point of view of notation, because now the even subsequence and the odd subsequence coincide with the even-indexed terms and the odd-indexed terms, respectively. 

There is an interesting relationship between the subsequences of type-1 and type-2 sequences, which had already been noticed in the case of Fibonacci numbers \cite{Raj04, Sil06}:   

\begin{proposition}	
	\label{prop:subsequences}
	Let $\big\{ G_n \big\}_{n=1}^{\infty}$ be a type-1-sequence and $\big\{ J_n \big\}_{n=1}^{\infty}$ a type-2-sequence. Then, 
	\begin{enumerate}
	\item The odd and even subsequences of $\big\{ G_n \big\}_{n=1}^{\infty}$ are of type 2.
	\item The odd and even subsequences of $\big\{ J_n \big\}_{n=1}^{\infty}$ are also of type 2.
	\end{enumerate}
\end{proposition}

\textbf{Proof:} In order to obtain a recurrence equation for the odd and even subsequences we adapt a simple (yet clever) technique developed in \cite{Sil06} for the Fibonacci numbers.  
\begin{enumerate} 
\item Let's start with the subsequences of $\big\{ G_n \big\}_{n=1}^{\infty}$. We have 
\begin{align}
	q G_n &= 2q G_n - q G_n \label{eq:part1A} \\
	p G_{n+1} &= p^2 G_n + pq G_{n-1}. \label{eq:part2A}
\end{align}
Adding (\ref{eq:part1A}) and (\ref{eq:part2A}) we get 
\begin{align}
	p G_{n+1} + q G_n &= (p^2+2q) G_n + pq G_{n-1} - q G_n \nonumber \\
	            G_{n+2} &= (p^2+2q) G_n + q(p G_{n-1} - G_n) \nonumber \\
	            G_{n+2} &= (p^2+2q) G_n - q^2 G_{n-2}. \label{eq:part3A} 
\end{align}
The recurrence equation \ref{eq:part3A} applies to both the odd and even subsequences of $\big\{ G_n \big\}_{n=1}^{\infty}$, and has the form of Equation \ref{eq:recurrence-type2}. To complete the proof of this case we just have to verify that $q \leq p$ effectively implies $q^2 < p^2+2q - 1$. 
\item Now we do the same for the subsequences of $\big\{ J_n \big\}_{n=1}^{\infty}$. We have 
\begin{align}
	q J_n &= 2q J_n - q J_n \label{eq:part1B} \\
	p J_{n+1} &= p^2 J_n - pq J_{n-1}. \label{eq:part2B}
\end{align}
Subtracting (\ref{eq:part1B}) from (\ref{eq:part2B}) we get 
\begin{align}
	p J_{n+1} - q J_n &= (p^2-2q) J_n + q J_n - pq J_{n-1} \nonumber \\
	           J_{n+2} &= (p^2-2q) J_n + q(J_n - p J_{n-1}) \nonumber \\
	           J_{n+2} &= (p^2-2q) J_n - q^2 J_{n-2}. \label{eq:part3B}
\end{align}
Again, note that Equation \ref{eq:part3B} has the form of Equation \ref{eq:recurrence-type2}. To complete the proof we just have to verify that $q^2 < p^2+2q - 1$ is equivalent $q < p-1$ when $q$ and $p-1$ are both positive. 
\end{enumerate}
\EndProof

From this we can derive the following

\begin{corollary}	
	\label{coro:odd-subsequences}
	Let $\big\{ G_n \big\}_{n=1}^{\infty}$ be a totally greedy type-1-sequence and $\big\{ J_n \big\}_{n=1}^{\infty}$ a totally greedy type-2-sequence. Then, the odd subsequences of $\big\{ G_n \big\}_{n=1}^{\infty}$ and $\big\{ J_n \big\}_{n=1}^{\infty}$ are also totally greedy. 
\end{corollary}

\textbf{Proof:} Proposition \ref{prop:subsequences} tells us that odd subsequences of type-1 or type-2 sequences are of type 2. On the one hand, the odd subsequences of $\big\{ G_n \big\}_{n=1}^{\infty}$ have the form 
\begin{displaymath}
	G'_k = G_{2k-1} = 
	\begin{cases}
		1 \mbox{ if } k=1, \\
		a' = pa+q \mbox{ if } k=2, \\
		p_1'G'_{k-1} - q'G'_{k-2} \mbox{ if } k>2,  
	\end{cases}
\end{displaymath}
where $p_1' = p^2+2q$ and $q' = q^2$. From Theorem \ref{theo:type2} we know that type-2 sequences are totally greedy for any $a$, and in particular for $a' = pa+q$, hence the result follows. 

On the other hand, the odd subsequences of $\big\{ J_n \big\}_{n=1}^{\infty}$ have the form 
\begin{displaymath}
	J'_k = J_{2k-1} = 
	\begin{cases}
		1 \mbox{ if } k=1, \\
		a' = pa-q \mbox{ if } k=2, \\
		p_2'J'_{k-1} - q'J'_{k-2} \mbox{ if } k>2,  
	\end{cases}
\end{displaymath}
where $p_2' = p^2-2q$ and $q' = q^2$. Again, the total greediness of $\big\{ J'_k \big\}_{k=1}^{\infty}$ follows directly from Theorem \ref{theo:type2}. 

\EndProof 

It is only at this point, after seeing the proof of Corollary \ref{coro:odd-subsequences}, that we can fully assess the need of introducing the parameter $a$ in our sequences, which has been the source of so much trouble throughout Sections \ref{sec:type1} and \ref{sec:type2}. 

Anyway, the even subsequences of $\big\{ G_n \big\}_{n=1}^{\infty}$ and $\big\{ J_n \big\}_{n=1}^{\infty}$ are more complicated because we have to insert the term $1$ at the beginning of the sequence, and this sort of \lq destroys\rq \ the recurrence. Nonetheless, we can still make some progress in some special cases, which are outlined below. 

To begin with, let $\big\{ G_n \big\}_{n=1}^{\infty}$ be a  type-1-sequence, and let $\big\{ G''_k \big\}_{k=1}^{\infty}$ denote the subsequence of its even terms, modified as follows: 
	\begin{equation}
	\label{eq:even-terms-type1-modified}
	G''_k =  
	    \begin{cases}
		    1 \mbox{ if } k=1, \\
		    a \mbox{ if } k=2, \\
		    p_1''G''_{k-1} - q''G''_{k-2} = G_{2k-2} \mbox{ if } k>2,  
	    \end{cases}
    \end{equation}
where $p_1'' = p^2+2q$ and $q'' = q^2$.  

\begin{corollary}	
	\label{coro:even-subsequences-of-type1}
	Let $\big\{ G_n \big\}_{n=1}^{\infty}$ be a totally greedy type-1 sequence with $p=q+1$ and $a=q$, and let $\big\{ G''_k \big\}_{k=1}^{\infty}$ be the subsequence of its even terms, as defined in Equation \ref{eq:even-terms-type1-modified}. Then $\big\{ G''_k \big\}_{k=1}^{\infty}$ is also totally greedy. 
\end{corollary}

\textbf{Proof:} We assume that $p$ and $q$ are fixed. From the equality $p_1''G''_{k-1} - q''G''_{k-2} = G_{2k-2}$ we get $\displaystyle a = \frac{pq+q^2}{2q-1}$, which clearly complies with the conditions of Theorem \ref{theo:type1}, i.e. $2 \leq a \leq p+q$. Thus, any combination of values of $p$ and $q$ such that $\displaystyle a = \frac{pq+q^2}{2q-1}$ is an integer will fit our purposes. In particular, $p=q+1$ and $a=q$ will do. 

\EndProof 

Note that the proof of Corollary \ref{coro:even-subsequences-of-type1} suggests how to look for other combinations of values for $p,q,a$. Another possibility would be to define Equation \ref{eq:recurrence-type1}, so that it starts with three given initial values, i.e. $1, a, b$, but that falls without the scope of this paper. 

Now let $\big\{ J_n \big\}_{n=1}^{\infty}$ be a  type-2-sequence, and let $\big\{ J''_k \big\}_{k=1}^{\infty}$ denote the subsequence of its even terms, modified as follows: 
	\begin{equation}
	\label{eq:even-terms-type2-modified}
	J''_k =  
	    \begin{cases}
		    1 \mbox{ if } k=1, \\
		    a \mbox{ if } k=2, \\
		    p_2''J''_{k-1} - q''J''_{k-2} = J_{2k-2} \mbox{ if } k>2,  
	    \end{cases}
    \end{equation}
where $p_2'' = p^2-2q$ and $q'' = q^2$.

\begin{corollary}	
	\label{coro:even-subsequences-of-type2}
	Let $\big\{ J_n \big\}_{n=1}^{\infty}$ be a totally greedy type-2 sequence with $a=p-q$, and let $\big\{ J''_k \big\}_{k=1}^{\infty}$ be the subsequence of its even terms, as defined in Equation \ref{eq:even-terms-type2-modified}. Then $\big\{ J''_k \big\}_{k=1}^{\infty}$ is also totally greedy. 
\end{corollary}

\textbf{Proof:} Again, we are assuming that $p$ and $q$ are fixed, and we are looking for a suitable value of $a$. From the equality $p_2''J''_{k-1} - q''J''_{k-2} = J_{2k-2}$ we get the desired condition $a=p-q$. 

\EndProof



\section*{Disclaimer and acknowledgements}
\label{sec:disclaimer}

This research did not receive any specific grant from funding agencies in the public, commercial, or
not-for-profit sectors. 



\end{document}